\documentclass[11pt, oneside]{article}   	% use "amsart" instead of "article" for AMSLaTeX format
\usepackage{geometry}                		% See geometry.pdf to learn the layout options. There are lots.
\usepackage{graphicx,color}				% Use pdf, png, jpg, or eps§ with pdflatex; use eps in DVI mode
\usepackage[bf,hang]{caption}
\usepackage{subcaption,url}								% TeX will automatically convert eps --> pdf in pdflatex		
\usepackage{amssymb, amsmath, amsthm}
\usepackage{authblk}

\newtheorem{theor}{Theorem}[section]

 \newtheorem{prop}{Proposition}[section]
 
\def\R{\mathbb{R}}

\def\N{\mathbb{N}}
\def\o{\overline{u}}

\allowdisplaybreaks

\begin{document}
\title{On the time fractional heat equation with obstacle}
\author[1]{C. Alberini}
\author[2]{R. Capitanelli}
\author[3]{M. D'Ovidio}
\author[4]{S. Finzi Vita}
%\date{}							% Activate to distplay a given date or no date

\affil[1]{{\it\footnotesize{Dipartimento SBAI, Sapienza Universit\`a di Roma, carlo.alberini@uniroma1.it}}}
\affil[2]{{\it\footnotesize{Dipartimento SBAI, Sapienza Universit\`a di Roma, raffaela.capitanelli@uniroma1.it}}}
\affil[3]{{\it\footnotesize{Dipartimento SBAI, Sapienza Universit\`a di Roma, mirko.dovidio@uniroma1.it}}}
\affil[4]{{\it\footnotesize{Dipartimento di Matematica, Sapienza Universit\`a di Roma, stefano.finzivita@uniroma1.it}}}

\date{\today}							% Activate to distplay a given date or no date

\maketitle
%\section{}
%\subsection{}

{\bf Abstract.}
%Si vuole studiare il comportamento delle soluzioni del problema parabolico dell'ostacolo in un dominio limitato con derivata temporale frazionaria (nel senso di Caputo) al variare del parametro $\alpha\in(0,1)$. Da un punto di vista numerico nel caso 1D si adattano combinandoli gli schemi alle differenze finite gi\`a testati per la {\it fractional heat equation} sulla base degli articoli di \cite{giga} e\cite{JLZ}, e per il problema parabolico classico con ostacolo, nella versione con Heaviside e in quella di \cite{BS}. 
We study a Caputo time fractional degenerate diffusion equation which we prove to be equivalent to the  fractional parabolic obstacle problem, showing that its solution evolves for any $\alpha\in(0,1)$ to the same stationary state, the solution of the classic elliptic obstacle problem. The only thing which changes with $\alpha$ is the convergence speed. \\
We also study the problem from the numerical point of view, comparing some finite different approaches, and showing the results of some tests.
These results extend what recently proved in \cite{ACF} for the case $\alpha=1$.
 
%\section{Introduction (da scrivere, precisando modello di partenza e propriet\`a)}

\section{Introduction}
%Il problema continuo da cui si parte \`e il seguente
%\begin{eqnarray}
%\label{model}
%\left\{\begin{array}{ll}
%\partial^\alpha_t u -H(u-\psi)\Delta u =0 \ &\hbox{ in } \Omega\times  (0,T) \\
%u(x,0)=u^0(x) &\hbox{ in } \Omega \\ 
%u=0 &\hbox{ on } \partial\Omega\times  (0,T), 
%\end{array}\right.
%\end{eqnarray}
%dove $\Omega$ \`e un dominio limitato di $\R^2$,  $u^0$ \`e il dato iniziale (che soddisfa le condizioni di Dirichlet al bordo), $\psi$ la funzione ostacolo (che verifica la condizione di compatibilit\`a $\psi\le 0$ al bordo), $H$ la funzione segno ($H(x)=1$ se $x>0$, $H(x)=0$ altrimenti), e $\alpha\in(0,1)$. La derivata temporale frazionaria nel senso di Caputo \`e definita come:
%\begin{equation}
%\label{frac_heat}
%\partial^\alpha_t u(x,t):=\frac{1}{\Gamma(1-\alpha)}\int_0^t (t-s)^{-\alpha}\partial_su(x,s)\ ds,
%\end{equation}
%dove $\Gamma$ \`e la funzione Gamma. 
In the works of the last decades, the use of fractional calculus has known an ever increasing use in describing phenomena from the most disparate scientific fields, from biology to mechanics, to superslow diffusion in porous media till financial interests.

Aim of this paper is to study the following problem, which generalizes the one addressed in \cite{ACF}, for $0 < \alpha <1$:
\begin{eqnarray}
\label{model1}
	\left\{
	\arraycolsep=1.4pt\def\arraystretch{1.5}
	\begin{array}{ll}
		\partial^\alpha_t u   -H\left(u-\psi\right)\Delta u =0 \ &\text{ a.e. in } \Omega, \text{ for all }t \in (0,T),	\\
		u\left(0\right) = u^{0} 						& \text{ in } \Omega,	\\
		u=0 & \text{ on } \partial\Omega, \text{ for all }t \in (0,T),
	\end{array}\right.
\end{eqnarray}
with $T>0$, where $H$ is the extended Heaviside function such that $H(0)=0$, that is

\begin{equation}
\label{heavi}
	H(r)=\left\{
	\arraycolsep=1.4pt\def\arraystretch{1.5}
	\begin{array}{ll} 1   &\text{ for } r>0\cr
		%[0,1]   &\text{ for } r=0 \cr 
		0  &\text{ for } r\le0
			\end{array}\right. 
\end{equation}
and $ \Omega$ is a bounded domain in $\R^n$,  $n\in\N$,  with smooth boundary. For simplicity, we omitted the presence of a forcing term in the problem.

Here $\partial^\alpha_t u$ denotes the Caputo fractional derivative, that is,
\begin{equation}
\label{frac_heat}
\partial^\alpha_t u(x,t):=\frac{1}{\Gamma(1-\alpha)}\int_0^t (t-s)^{-\alpha}\partial_su(x,s)\ ds,
\end{equation}
where  $\Gamma$ is the  Gamma function. 

The Riemann-Liouville derivative is defined as
\begin{align}\label{RL}
^R\partial^\alpha_t u(x,t) := \frac{\partial}{\partial t}\frac{1}{\Gamma(1-\alpha)}\int_0^t (t-s)^{-\alpha}\, u(x,s)\ ds .
\end{align}
We observe that
\begin{align*}
\, ^R\partial^\alpha_t C = \frac{1}{\Gamma(1-\alpha)} \partial_t \int_0^t \frac{C}{(t-s)^\alpha} ds = C \frac{t^{-\alpha}}{\Gamma(1-\alpha)}
\end{align*}
and the following relation between derivatives holds
\begin{align}
\label{RelDerProp22}
\partial^\alpha_t u(x,t) =\, ^R\partial^\alpha_t \big( u(x,t) - u(x,0)\big) = \, ^R\partial^\alpha_t u(x,t) - u(x,0) \frac{t^{-\alpha}}{\Gamma(1-\alpha)}.
\end{align}
We refer to the book \cite{Dieth} for further details on fractional derivatives.

Let us consider the complementarity system, for all $t>0$,
\begin{eqnarray}
\label{parab_cs}
	\left\{
	\arraycolsep=1.4pt\def\arraystretch{1.5}
		\begin{array}{ll}
		w(x,t) \geq \psi(x) 								&\  	\text{in } \Omega		\\
		\partial^\alpha_t w \geq \Delta w 							&\ 	\text{a.e. in } \Omega		\\
		\displaystyle (w-\psi)(\partial^\alpha_t w-\Delta w ) = 0			&\	\text{a.e. in }\Omega		\\
		w\left(x,0\right)=u^0 							&\	 \text{in } \Omega		\\
		w(x,t) = 0	&\	 \text{on } \partial\Omega .
		\end{array}
	\right. 
\end{eqnarray}

In Section 2 we will prove that
\begin{description}
	\item[-] by analogy with the classical case, (\ref{model1}) is equivalent to the complementarity system \eqref{parab_cs};
	\item[-] asymptotically, the solution evolves for each $\alpha$ towards the same steady state of the classical problem, that is
		\begin{equation}
			\label{BS}
			\bar{u}\ge\psi,\quad - \Delta \bar{u} \ge 0,\quad (\bar{u} - \psi) \Delta \bar{u} = 0 \quad\hbox{in }\Omega\times (0, T)
		\end{equation}
with different convergence speed (for $\alpha = 1$ exponential, for $\alpha\in (0,1)$ polynomial).
\end{description}

In Section 3 we present three possible finite difference schemes for the numerical approximation of problem (\ref{model1}) or its equivalent form given by the complementary system (\ref{parab_cs}). The Caputo derivative has been discretized using standard methods found in the literature, the so-called L1 or Convolution Quadrature (CQ) approaches (see e. g. \cite{JLZ} or \cite{giga}). The space discretization has been carried out through the semi implicit finite differences scheme introduced in \cite{ACF} for problem (\ref{model1}), or through the implicit scheme of \cite{BS} for the evolutive obstacle problem (\ref{parab_cs}). Note that for $\alpha\to 1^{-}$ all these schemes give back the known results of the classic heat equation with obstacle.

Our aim was not to find an optimal strategy of approximation for the problem, but only to derive working schemes in orders  to confirm through explicit simulations the behavior of the solution as  characterized by the results of Section 2. For that, in Section 4, we have tested the previous schemes on a couple of one-dimensional examples. Such simulations also allow to compare their reliability and  computational cost. The semi implicit approach requires time step restrictions (strongly increasing as $\alpha \to 0^+$) in order to detect the correct contact set with the obstacle, restrictions unnecessary for the implicit approach.
On the other hand, in the first case each time iteration is much less expensive, since it reduces to a single linear system solution. So, after all, all the proposed schemes appear to be competitive.

\medskip For the following, we assume that the  initial datum $u^0$ and  the (independent of time) obstacle $\psi$ satisfy the following conditions
\begin{equation}
\label{obst}
	 u^0\in H^1_0\left(\Omega\right),\, \psi\in H^2\left(\Omega\right),\ \psi \le 0 \text{ on } \partial\Omega.
\end{equation}

We define as solution of problem (\ref{model1}) a function $u\in L^{2}\left(0,T;H^1_0\left(\Omega\right)\cap H^2\left(\Omega\right)\right)$ with $\partial _t^\alpha u(x,t)\in L^{2}\left(0,T;L^2(\Omega)\right)$ which solves problem (\ref{model1}).

\section{Equivalence of the problems and asymptotic behaviour}

In the present section we analyze problem (\ref{model1}), showing that, under suitable conditions,   it is equivalent to  complementarity system \eqref{parab_cs}. Let us introduce the following hypoteses:
\begin{description}
 	\item[H$_{1}$:] $u^{0} > \psi$ a.e. in $\Omega$;
	\item[H$_{2}$:] $\Delta \psi  \leq 0$ a.e. in $\Omega$.
\end{description}

\begin{prop}\label{equiv} 
Assume that conditions \eqref{obst}, {\bf H}$_{1}$ and {\bf H}$_{2}$ hold. Then  problems \eqref{model1} and \eqref{parab_cs} are equivalent.

\end{prop}

\noindent{\bf Proof.}

First we note that  if $w$ solves problem    \eqref{parab_cs}, then $w$ solves \eqref{model1}.   
In fact, initial and boundary conditions are the same in \eqref{model1} and \eqref{parab_cs}. Since $(w-\psi)(\partial^\alpha_t w-\Delta w ) = 0$ a.e. in $\Omega$,  
where $w>\psi$, then we obtain $\partial^\alpha_t w -\Delta w = 0$  a.e $\text{in } \Omega$ whereas, where $w(x,t)=\psi(x)$ then $\partial_t w=0$ which implies $\partial_t^\alpha w=0$. 

Now we prove that if $u$ solves problem \eqref{model1} then $u$ coincides with  solution $w$ of \eqref{parab_cs}.
Initial and boundary conditions are the same in \eqref{model1} and \eqref{parab_cs}. 
Let us now prove that if  $u(x,t)$ is a solution of \eqref{model1} then necessarily $u(.,t)\geq \psi(.)$ in $\Omega$ for any time $t$. Assume that $u< \psi$. Then $\partial^\alpha_t u=0$ from \eqref{model1}. From ${\bf H}_1$, $u^0>\psi$. Moreover, due to the fact that $u < \psi < u^0$, we obtain from \eqref{RelDerProp22}
\begin{align*}
\partial^\alpha_t u = \, ^R\partial^\alpha_t ( u - u^0 ) < 0,
\end{align*}
which does not agree with $\partial^\alpha_t u=0$. Then $u\geq \psi$ by contradiction and the first inequality of (\ref{parab_cs}) holds.
The  equation in the third line of \eqref{parab_cs} is trivially satisfied where $u(x,t)= \psi(x)$; where $u(x,t)>\psi(x),$ from (\ref{model1})  we obtain $\partial^\alpha_t   -\Delta u =0$, so it is always true.
Concerning the second inequality of \eqref{parab_cs}, we have already seen that it is satisfied (with the equal sign) when $u> \psi$. But when $u(x,t)=\psi(x)$,  (\ref{model1}) and  assumption {\bf H}$_{2}$ imply that
$$ \partial^\alpha_t  u  -\Delta u = -\Delta \psi  \geq 0 .$$
Then $u$ solves problem \eqref{parab_cs}.\qed

\bigskip

We recall that the same problem for $\alpha=1$ has been faced in \cite{ACF}, proving in
 particular   the  equivalence of  problem (\ref{model1}) with a parabolic obstacle problem.
Let $v$ be the unique solution  of  problem \eqref{parab_cs} with $\alpha=1$, that is
\begin{eqnarray}
\label{parab_cs1}
	\left\{
	\arraycolsep=1.4pt\def\arraystretch{1.5}
		\begin{array}{ll}
		v(x,t) \geq \psi(x) 								&\  	\text{in } \Omega		\\
		 \partial_t v \geq \Delta v 							&\ 	\text{a.e. in } \Omega	\\
		\displaystyle (v-\psi)( \partial_t v -\Delta v ) = 0			&\	\text{a.e. in }\Omega		\\
		v\left(x,0\right)=u^0 								&\	 \text{in } \Omega		\\
		v(x,t) = 0.										&\	 \text{on } \partial\Omega
		\end{array}
	\right. 
\end{eqnarray}

Let us introduce
\begin{equation}
\label{sol}
u(x,t) = \int_0^\infty v(x,s)\, \ell(s,t)\, ds
\end{equation}
where the probability density $\ell : [0, \infty) \times [0, \infty) \to [0, \infty)$ satisfies (see \cite{Dov})
\begin{equation*}
\begin{cases}
\displaystyle \partial^\alpha_t \ell = - \frac{\partial \ell}{\partial s}, \quad (s,t) \in (0,\infty) \times (0, \infty)\\
\displaystyle \ell(s,0) = \delta(s), \quad s>0,\\
%&\displaystyle \ell(0,t) = \frac{t^{-\alpha}}{\Gamma(1-\alpha)}, \quad t>0\\
\displaystyle \ell(s,t) =0, \quad s < 0.
\end{cases}
\quad \textrm{and} \quad
\begin{cases}
\displaystyle ^R\partial^\alpha_t \ell = - \frac{\partial \ell}{\partial s}, \quad (s,t) \in (0, \infty) \times (0, \infty),\\
\displaystyle \ell(s,0) = \delta(s), \quad s>0,\\
\displaystyle \ell(0,t) = \frac{t^{-\alpha}}{\Gamma(1 - \alpha)}, \quad t>0,\\
\displaystyle \ell(s,t) =0, \quad s<0.
\end{cases}
\end{equation*}

\newpage

\begin{prop}\label{equiv2} 
The function $u$ defined in \eqref{sol} solves \eqref{parab_cs}.
\end{prop}

\noindent{\bf Proof.}

 As $v \ge \psi$, then
\begin{align*}
\int_0^\infty v(x,s)\, \ell(s,t)\, ds \geq \int_0^\infty \psi(x)\, \ell(s,t)\, ds = \psi(x),
\end{align*}
and we obtain that $u\ge \psi.$

Now we prove that $(u-\psi)(\partial^\alpha_t u-\Delta u)=0.$ 
For $u> \psi$, it holds the classical theory about  fractional Cauchy problems (see \cite{Baz2000}, \cite{Koc89} and \cite[Theorem 5.2]{CapDov}) and therefore we have that, in $L^2(\Omega)$,
\begin{align*}
\partial_t^\alpha u(x,t) = \Delta u(x,t), \quad u(x,0)= u^0(x).
\end{align*}

Concerning the second inequality of (\ref{parab_cs}), we observe that
 if $\partial_t v - \Delta v \ge 0$, then
\begin{align*}
\Delta \int_0^\infty v(x,s) \, \ell(s,t)\,  ds  
\le & \int_0^\infty \partial_s v(x,s) \, \ell(s,t)\,  ds\\
= & - v(x,0) \frac{t^{-\alpha}}{\Gamma(1-\alpha)} - \int_0^\infty v(x,s)\, \partial_s \ell(s, t)\, ds\\
= & - v(x,0) \frac{t^{-\alpha}}{\Gamma(1-\alpha)} + \, ^R\partial^\alpha_t u(x,t).
\end{align*}
Since $u(x,0)=v(x,0)$ from the relation \eqref{RelDerProp22} we obtain  
$\partial^\alpha_t u - \Delta u \ge 0$.

The initial condition and the boundary condition for the function $u$ follow trivially from  the initial condition and the boundary condition for the function $v$.
\qed

\medskip

In \cite{ACF} the asymptotic solution of  \eqref{parab_cs1} has been characterized as the solution of the corresponding stationary (elliptic) obstacle problem.
In particular,   it has been proved that 
under conditions   \eqref{obst}, {\bf H}$_{1}$ and {\bf H}$_{2}, $  the  solution  $v(t)$  to the  parabolic  obstacle problem  \eqref{parab_cs1} converges strongly in $H^1_0(\Omega),$
for $t\to\infty,$ to the unique solution  $\overline{u}$  of the corresponding  stationary obstacle problem 
\begin{equation}
	\label{ell_obst}
	\overline{u}\in H^1_0(\Omega),\quad \overline{u}\geq \psi,\quad -\Delta\overline{u}\geq 0,\quad (\overline{u}-\psi)(\Delta \overline u )=0\ \hbox{ a.e. in }\Omega.
\end{equation}

Moreover, it has been proved that  there is a constant $C > 0$ such that, for every $t \geq 1$,
\begin{equation}
\label{stima1}
	\left\Vert v(t) - \overline{u}\right\Vert_{H^1\left(\Omega\right)} \leq e^{-Ct} \ .
\end{equation}

In the next Theorem we prove that a similar  result holds for any $\alpha\in(0,1).$

\begin{theor}\label{main} 
Assume conditions   \eqref{obst}, {\bf H}$_{1}$ and {\bf H}$_{2}$ hold. Let $u(t)$ be  the function defined in \eqref{sol} solution of  \eqref{parab_cs}.
Let $\overline{u}$  be the unique solution of the obstacle problem  \eqref{ell_obst}. 
Then, $u(t)$  converges to $\overline{u}$ 
%strongly in $H^1_0(\Omega)$
for $t\to\infty,$ and there is a constant $C > 0$ such that, for every $t \geq 1$,
\begin{equation}
\label{stima1}
	\left\Vert u(t) - \overline{u}\right\Vert_{L^1\left(\Omega\right)} \leq \kappa_\Omega\,  E_\alpha(- C t^\alpha) \ 
\end{equation}
for any $\alpha \in (0,1),$ where $\displaystyle\kappa_{\Omega} = \left( \int_\Omega dx \right)^{1/2} $ and
\begin{align}
\frac{E_\alpha(-C t^\alpha)}{J(t)} \to 1, \quad \textrm{as} \quad  t\to \infty,
\end{align} 
with
\begin{align}\label{jei}
J(t) = \frac{1}{C} \frac{ t^{-\alpha}}{\Gamma(1-\alpha)}.
\end{align}

\end{theor}

\noindent{\bf Proof.}
Since 
\begin{align*}
\| v - \bar{u} \|_{L^1(\Omega)} \leq  \kappa_{\Omega} \; \| v - \bar{u} \|_{L^2(\Omega)} 
\end{align*}
we have 
\begin{align*}
\| u(t) - \bar{u} \|_{L^1(\Omega)
} 
= & \int_\Omega \bigg| \int_0^\infty v(x,s)\, \ell(s,t)\, ds  - \bar{u}(x)\bigg| dx \\
= & \int_\Omega \bigg| \int_0^\infty \left( v(x,s) - \bar{u}(x)\right) \ell(s,t)\, ds \bigg| dx\\
\le & \int_\Omega  \int_0^\infty \bigg| v(x,s) - \bar{u}(x)\bigg| \ell(s,t)\, ds  dx\\
= & \int_0^\infty \| v(x,s) - \bar{u}(x) \|_{L^1(\Omega)} \, \ell(s,t)\, ds\\
\le & \kappa_\Omega \int_0^\infty e^{-C s} \, \ell(s,t)\, ds\\
= & \kappa_\Omega \, E_\alpha(- C t^\alpha), \quad \forall\, \alpha \in (0,1)
\end{align*}
where, in the last step, we used the fact that $E_\alpha$ is the Laplace transform of the density $\ell$. For the asymptotic behaviour of the Mittag-Leffler, consult the book \cite[formula (4.4.17)]{GKMR-book}.

\qed

\section{Numerical approximation }

%%%%%%%
%The aim of this section is not to find an optimal strategy of approximation for our problem, but only to derive a working scheme in orders  to show through explicit simulations the behavior of the solution as already characterized by the results of the previous section.

The starting idea for a numerical approximation of the problem  has been to combine  classical  time-stepping discretization schemes for the Caputo derivative, such as the convolution quadrature (CQ) or the finite difference (L1) scheme, see e.g. \cite{JLZ}, with schemes usually working for the parabolic obstacle problem (\ref{parab_cs}), such as the one proposed in \cite{BS}, but also the semi-implicit f.d. scheme tested in \cite{ACF} for the equivalent Heaviside function formulation (\ref{model1}) of the problem. 

For sake of simplicity we work in the one-dimensional case, with  $\Omega=(a,b)$.
Let us call $h>0$  the space discretization step ($h=(b-a)/N$, so that we have $(N-1)$ internal nodes in $\Omega$, $x_i=a+ih$, for $i=1,...,N-1$) and $\tau=T/M$ the time discretization step (with $M$ time instants  $t^m=m\tau$, for $m=1,...,M$); with $\gamma_\alpha=\tau^\alpha/h^2$ we denote the  parabolic ratio between the steps related to a specific $\alpha$. Note that for fixed steps $h$ and $\tau$, if $\alpha$ decreases to zero  then $\gamma_\alpha$ quickly grows: in other words for small $\alpha$, in order to keep  $\gamma_\alpha$ small on a fixed mesh, the step $\tau$ has to be considerably reduced, with a significant increase of computational costs.

Here we analyze three possible approaches:

\begin{enumerate}
\item {\bf Scheme S1: solves problem (\ref{model1}) with L1 for the Caputo derivative and the semi implicit f.d. scheme of \cite{ACF} in space.}

The time discretization of the Caputo derivative by the L1 scheme leeds to the formula (see \cite{JLZ}):
$$
\partial^\alpha_t u(x,t^m) \simeq \frac{1}{\Gamma(2-\alpha)\tau^\alpha}\left\{u(x,t^m)-\sum_{k=0}^{m-1} C_{m,k}u(x,t^k)\right\}\ ,$$
where
$$ C_{m,0}:= f(m), \quad C_{m,k}:=f(m-k)-f(m-(k-1)) \ \ \hbox{ for } k=1,...,m-1 , $$
and
$$ f(r):= r^{1-\alpha}-(r-1)^{1-\alpha}, \hbox{ for } r\ge 1 .$$
Then, after  semidiscretization in time, we need to solve for any instant $t^m$ the equation
\begin{equation}
\label{m_eq}
\frac{1}{\Gamma(2-\alpha)\tau^\alpha}\left\{u(x,t^m)-\sum_{k=0}^{m-1} C_{m,k}u(x,t^k)\right\}-H(u-\psi)u_{xx} =0 \quad \hbox{ in } \Omega,
\end{equation}
with the same boundary conditions of (\ref{model1}). Since  $\Omega$ was splitted in $N$ subintervals through the nodes $x=\{x_i\}_i$, the initial data will be the vector  $u^0\in\R^{N-1}$, with $u^0_i=u^0(x_i)$;   applying the  semi implicit finite difference scheme of \cite{ACF} in space, the solution $u^1$ at the first discrete time $t^1=\tau$ will solve at any node the relation:
$$\frac 1 {g\tau^\alpha} \left(u^1_i-C_{1,0}u^0_i\right)=H(u^0_i-\psi_i)\delta^2 u^1_i:=H(u^0_i-\psi_i) \frac{u^1_{i-1}-2u^1_i+u^1_{i+1}}{h^2},$$
with $g=\Gamma(2-\alpha)$ (note that $0.8862\le g \le 1, \forall \alpha\in[0,1]$, so that this term will be negligible with respect to $\gamma_\alpha$). If we set $v^k_i=u^k_i-\psi_i$, redistributing all the  terms between the two members, it is equivalent to solve 
$$u^1_i -g\tau^\alpha H(v^0_i) \delta^2 u^1_i=C_{1,0}u^0_i \quad \hbox{ for every } i;$$
with vector notations it means that $u^1$ solves the linear system:
$$B^0u^1:=(I + g\gamma_\alpha H(v^0)*A)u^1=C_{1,0}u^0 ,$$
where $A$ is the usual tridiagonal matrix $(N-1)\times(N-1)$ with values $2$ on the main diagonal and  $-1$ on the two adjacent diagonals, and we denoted
$$ \{(H(v)*A)u\}_i:=H(v_i)(Au)_i=H(v_i)\sum_{j=1}^{N-1}A_{i,j}u_j\ .$$ 
Since the discrete solution could overstep the obstacle at some nodes, in particular when a large value of $\gamma_\alpha$ is used, the following correction is needed at any iteration:
$$u^1_i=\max(u^1_i,\psi(x_i))\ .$$
%si rende necessaria perch\'e in alcuni nodi la soluzione pu\`o superare l'ostacolo anche se di poco: ci\`o avviene usando l'$H$ esatta, o la sua approssimazione $\eta_n$ con $n$ molto grande. Da approfondire se esistano condizioni sui passi per cui lo schema non possa farlo.

In the same way we see that  $u^2$ solves for any $i$
$$\frac 1 {g\tau^\alpha} \left(u^2_i-C_{2,0}u^0_i-C_{2,1}u^1_i\right)=H(u^1_i-\psi_i)\delta^2 u^2_i,$$
that is the system
$$B^1u^2=(I + g\gamma_\alpha H(v^1)*A)u^2=C_{2,0}u^0+C_{2,1}u^1 ,$$
with the same matrix $A$ and the subsequent correction. In general, at any time step $u^m$ solves the linear system
\begin{equation}\label{sistema1}
B^{m-1} u^m=b^m,
\end{equation}
where we have set
\begin{equation}
\label{coefs1}
B^{m-1}=I + g\gamma_\alpha H(v^{m-1})*A, \qquad b^m=\sum_{k=0}^{m-1} C_{m,k}u^k\ ,
\end{equation}
followed by the correction
\begin{equation}
\label{correction}
u^m_i=\max(u^m_i,\psi(x_i)).
\end{equation}
Note that all the matrices $B^m$ are   symmetric positive definite (and M-matrices), since so it is  $A$, while  $H(v)\ge 0$. Then all the previous linear systems are well posed. 

With respect to the classic parabolic obstacle problem approach discussed in \cite{ACF}  there is  here an important difference. In that case (which corresponds to the case $\alpha=1$), when the solution at time $t^{m-1}$ touches the obstacle at node $x_i$, then $v^{m-1}_i=0$, $H(v^{m-1}_i)=0$, and (\ref{sistema1}) trivially yields:
$$ u^m_i=u^{m-1}_i \ ;$$
in other words, once touched the obstacle at a particular node the solution  does not change  anymore there, but only at the remaining free nodes. 
%Del resto  gi\`a nel problema continuo dove la $H$ si annulla non c'\`e pi\`u diffusione e $u_t=0$. 
In the general case of $\alpha\in(0,1)$ on the contrary,  at the contact time system (\ref{sistema1}) immediately yields for the $i$-th component:
$$u^m_i=\sum_{k=0}^{m-1} C_{m,k}u^k_i=C_{m,0}u^0_i+...+C_{m,m-1}\psi_i>\psi_i \sum_{k=0}^{m-1} C_{m,k}=\psi_i ,$$
since at least $u^0_i>\psi_i$ and $\sum_{k=0}^{m-1} C_{m,k}=1$. It follows that the solution has a little rebound at $x_i$ which detaches it again from the obstacle, and produces an (innatural) oscillating evolution from that time on. The width of such rebound depends on the size of $\gamma_\alpha$. Then, even if in principle the semi implicit scheme does not not require stability restrictions on the discretization steps, a small value of $\gamma_\alpha$ will be necessary to reduce the oscillations (they would vanish for $\tau\to 0$). %, since in the continuous setting a null Caputo derivative implies  $u_t=0$, then a constant solution in time). As a consequence of that,  the scheme would become enormously expensive, the more the more $\alpha$ is close to zero. 
The way to solve this difficulty is to remove the memory effect at the contact nodes of the solution with the obstacle, that is where  the obstacle retains the solution. This suggests to modify the scheme replacing the vector  $b^m$ of (\ref{coefs1})  by the vector $\hat{b}^m$ defined by
\begin{equation}
\label{coefs2}
\hat{b}^m=\max(H(v^{m-1})*b^m,u^{m-1});
\end{equation}
if $H(v^{m-1}_i)=0$, that is $u^{m-1}_i=\psi_i$, then $\hat{b}_i^m=u^{m-1}_i$ and $u^m_i=u^{m-1}_i$; otherwise $\hat{b}_i^m=b_i^m$ , since  $b_i^m\ge u^{m-1}_i$, and all remains as before. No rebound is still possible after a contact.

%Another possibility should be to adopt a variable time step approach for the scheme, at least after the first contact with the obstacle, not easy to formulate due to the memory effect of this model. In the present paper we will not consider this approach.  ({\rs togliere?})

\item  {\bf Scheme S2: solves problem (\ref{model1}) with CQ for the Caputo  derivative and the semi implicit f.d. scheme of \cite{ACF} in space.}

In this case the Caputo derivative is approximated through the so-called {\it convolution quadrature} (CQ) method, proposed by Lubich for the discretization of Volterra integral equations. In particular,  for a function $\phi(t)$ (with $\phi(0)=0$),  the Riemann-Liouville derivative $^R\partial^\alpha_t \phi$ defined in (\ref{RL}) can be approximated by the discrete convolution:
$$ ^R\partial^\alpha_\tau  \phi^m:= \frac{1}{\tau^\alpha}\sum_{j=0}^m c_j \phi^{m-j},$$
where $\phi^m=\phi(t_m)$, and the coefficients $\{c_j\}$ are obtained from a suitable power series expansion, connected to a specific approximation method for the ODE (see \cite{JLZ}). In the case of the Euler  backward method, it is known as the  Grunwald-Letnikov approximation, and provides the following recursive formula for the coefficients: 
$$c_0=1,\quad c_j=-\frac{\alpha-j+1}{j}c_{j-1}.$$
Then, using the  relation (\ref{RelDerProp22}) between the Caputo and the Riemann-Liouville derivatives we can rewrite the initial problem as 
$$ ^R\partial^\alpha_t (u-u^0)-H(u-\psi)\Delta u =0,$$
which discretized in time and space (with the same notations of S1) becomes:
$$\frac{1}{\tau^\alpha}\sum_{j=0}^m c_j (u^{m-j}-u^0)+\frac 1 {h^2}H(u^{m-1}-\psi)Au^m=0\ ,$$
equivalent to the solution at any iteration of the  linear system
\begin{equation}\label{sistema4}
B^{m-1} u^m=b^m,
\end{equation}
where this time we have set
\begin{equation}
\label{coefs4}
B^{m-1}=I + \gamma_\alpha H(u^{m-1}-\psi)*A, \qquad b^m=u^0-\sum_{j=1}^{m-1}c_j(u^{m-j}-u^0)\ ,
\end{equation}
followed again by the correction
\begin{equation}
\label{correction}
u^m_i=\max(u^m_i,\psi(x_i))\ .
\end{equation}
Even in this case the vector $b^m$ has to be modified inside the contact set in order to remove the memory effect and prevent rebounds, as done in (\ref{coefs2}). In fact when $u^{m-1}_i=\psi_i$ then again $b^m_i=0$, and from (\ref{coefs4}) we get
$$u^m_i=u^0_i-\sum_{j=1}^{m-1}c_j(u^{m-j}_i-u^0_i)>u^0_i-(\psi_i-u^0_i)\sum_{j=1}^{m-1}c_j>\psi_i\ ,$$
since $u^{m-j}_i\ge\psi_i$ for any $j$, and $\sum_{j=1}^{m-1}c_j\ge-1$.
% ({\bf controllare})

\item {\bf Scheme S3: solves problem (\ref{parab_cs}) with L1 for the Caputo derivative and the scheme of \cite{BS} for the evolutive obstacle problem.}

If we discretize the equation of system (\ref{parab_cs}) through  finite differences, using the L1 scheme for the Caputo derivative,  we get the equation
$$(u^m-\psi)^T\left(\frac 1 {g\tau^\alpha}(u^m-\sum_{k=0}^{m-1}C_{m,k}u^k)+\frac 1 {h^2}Au^m\right)=0\ .$$
Setting $y^m=u^m-\psi$, and remembering that $\sum_{k=0}^{m-1}C_{m,k}=1$, it is equivalent to
$$y^m \left(y^m+ g\gamma_\alpha A(y^m+\psi) -\sum_{k=0}^{m-1}C_{m,k}y^k\right)=0 \ .$$
Then $y^m=\max(0,x^m)$ is solution of the previous equation if $x^m$ solves
\begin{equation}
\label{picardC}
(I+ g\gamma_\alpha AP(x))x=b^m\ ,
\end{equation}
where now
\begin{equation}
\label{coefs3}
b^m=\sum_{k=0}^{m-1}C_{m,k}y^k-g\gamma_\alpha A \psi=\sum_{k=0}^{m-1}C_{m,k}u^k-\psi- g\gamma_\alpha A \psi\ ,
\end{equation}
while $P(x)=diag(p(x_i))$ denotes the diagonal matrix with  $p(x_i)=1$ if $x_i>0$ and $p(x_i)=0$ otherwise.
As seen in \cite{BS},  (\ref{picardC}) can be solved by  the so-called Picard iterations:
$$P^0=O, \quad (I+g\gamma_\alpha AP^n)x^{n+1}=b^m,\quad P^{n+1}=diag(p(x^{n+1}) )\quad \hbox{ per } n=0,1,... $$
until $P^{n+1}=P^n$ ($O$  is the null matrix); at that point $x^{n+1}$ is the  sought solution.

\end{enumerate}

\medskip Of course other schemes could be obtained by different combinations of specific numerical approaches, but for our purposes the three previous schemes were sufficient to perform in the next section explicit simulations of the problem. 

\section{Numerical tests}

We have applied the schemes described in Section 3 to some specific examples, for different values of $\alpha$. We have choosen a sufficiently large final time $T$, but also added a stopping time criterium in order to put in evidence the convergence towards the asymptotic solution. Since this convergence corresponds to the stabilization of the solution vector and to the satisfaction of the asymptotic complementarity relation 
$$ (u-\psi)\Delta u=0 ,$$
which  means $u$ harmonic (linear in 1D) outside the contact set, we have used the criterium
\begin{equation}\label{stop}
 STOP\quad  when\quad  \max(\|u^m-u^{m-1}\|_\infty,\|\left(u^m-\psi)A u^{m}\right\|_\infty) < tol\ ,
\end{equation}
for a given tolerance parameter $tol$. 

\begin{figure}[h!]
	\centering      	
	\includegraphics[width=7.5cm]{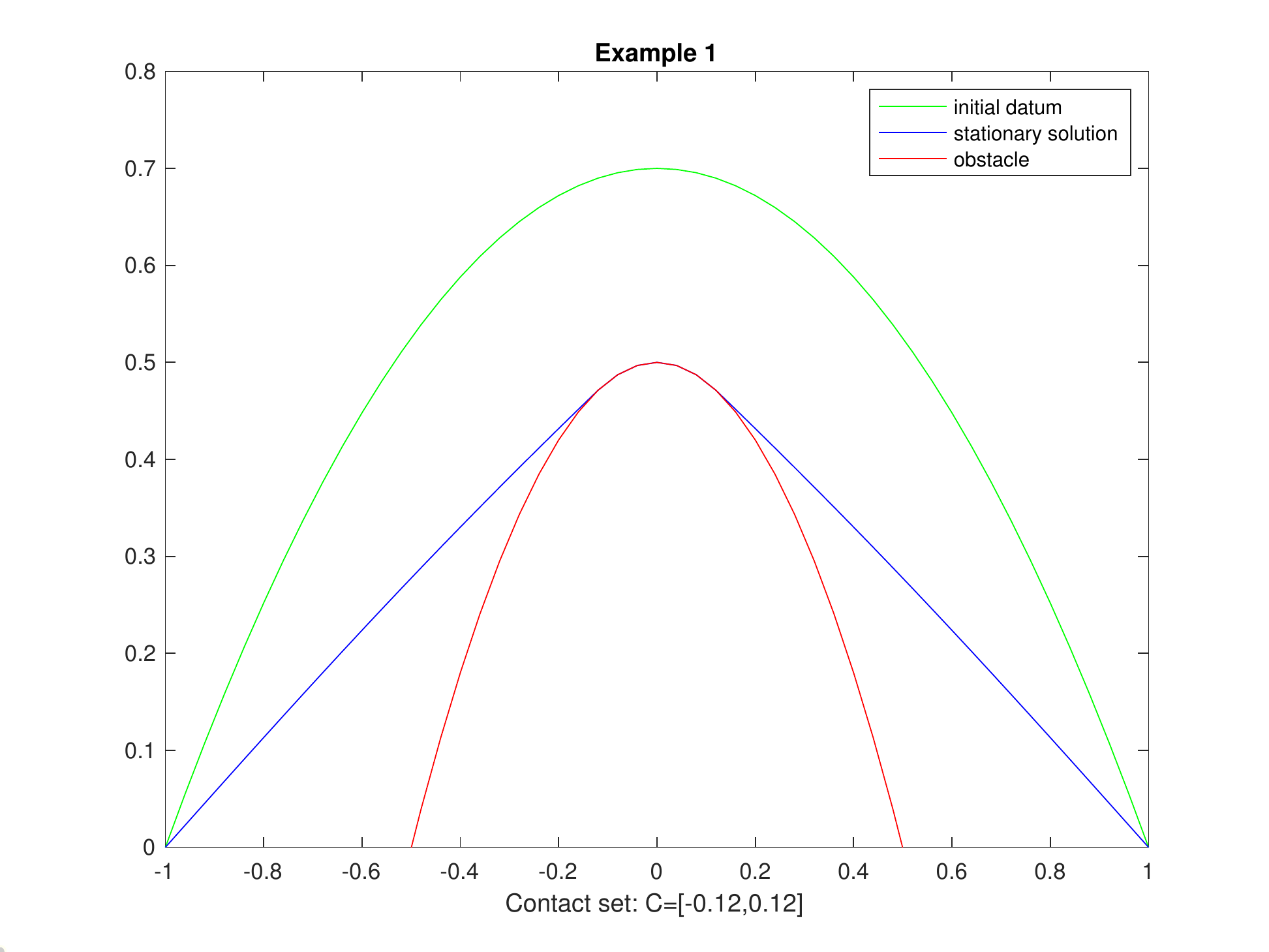}
	%\vskip-1.5cm
\caption{Data of Example 1.}
	\label{Ex1}
\end{figure}

\medskip\noindent {\bf Example 1.} $\Omega=(-1,1)$, $u^0(x)=0.7-0.7x^2$, $\psi(x)=0.5-2x^2$ (see Fig.\ref{Ex1}).

%\medskip Il programma ({\it heavifrac.m}) chiede di scegliere la costante $\alpha$,  il tempo finale $T$ (in questo esempio $T=1$ \`e sufficiente per osservare il contatto con l'ostacolo, mentre serve un $T$ molto pi\`u grande per avere la stabilizzazione della soluzione), e gli intervalli uniformi di tempo e spazio ($M$ e $N$), che determinano $\gamma_\alpha$. Chiede inoltre di scegliere quale schema applicare tra i precedenti. 
%Come eventuale criterio di arresto prima del tempo $T$ si era usato inizialmente quello della stabilizzazione delle iterate:
%$$ \| u^m - u^{m-1}\| < tol\ ,$$
%dove in genere si \`e considerato $tol=10^{-4}$. Si vede per\`o dai test che questo criterio non garantisce che la soluzione numerica finale sia sufficientemente vicina a quella stazionaria, che dovrebbe soddisfare la relazione di complementarit\`a (\ref{BS}): $ (u-\psi)\Delta u=0$, cio\`e $u$ armonica (lineare in 1D) fuori dall'insieme di contatto. Per questo motivo si \`e adottato il nuovo criterio, che viene soddisfatto in genere in tempi molto pi\`u lunghi:
%
%
%Lo schema S1 \`e stato testato, ma per quanto scritto nell'Osservazione 2, abbandonato. Le oscillazioni della soluzione a ridosso dell'ostacolo sono innaturali, e portano pure un altrettanto innaturale aumentato contatto della soluzione con l'ostacolo.

\medskip In the following Table~\ref{default1} we reported some results obtained by the simulations with schemes S1, S2 and S3, for different values of $\alpha$, $N$ and $\gamma_\alpha$, with $tol=10^{-4}$. We adopted the following notations:
\begin{itemize}
\item FC time = full contact time, that is the first time at which the solution has reached the whole contact set (no more changing in the successive iterations);
\item STOP time = the exit time according to  criterium (\ref{stop});
\item \# iter. = final number of time iterations;
\item \# Pic. = average number of  Picard iterations for each time step in scheme S3;
\item \# LS = approximate number of linear systems to be solved (essentially the product of the previous two quantities)
\item S =  working schemes (the ones detecting the correct contact set).
\end{itemize}
Looking at the table some remarks and comments are possible:

\begin{table}[ht]
\caption{}
\begin{center}
\begin{tabular}{|c|c|c|c|c|c|c|c|c|c|}
\hline
$\alpha$ & $N$ & $\gamma_\alpha$ & $\tau$ &FC time & STOP time & \# iter. & \# Pic. & \# LS& S \cr
\hline
   0 & 32 & 256 & any & first it. & no conv.  & 1 & 14 & & S3 \cr
       & 64 &1024 & any & first it. & no conv.  & 1 & 29 &  &S3\cr
      & 128 &4096 & any & first it. & second it. & 1 & 57 & &S3\cr
\hline
0.3 & 32 & 60& 0.0079 & 0.10 & $93.16$ & 11739 &12 & 140868 & all \cr
     & 32 & 75 & 0.016   & 0.11 & $93.15$ & 5580& 13 & 72540 & S3 \cr
     & 64 & 220 & 0.0059 & 0.49 & 1.33 & 226 & 21& 4746 & all \cr
     & 128 &850 & 0.005 & 0.2 & 1.66& 315 &  42& 13230 &all\cr
\hline
0.5 & 32 & 25 &0.009 & 0.16 & 8.38 & 880 & 8 & 7040& all\cr
   &  32 & 50 &0.038 & 0.19 & 8.39 & 221 & 11 & 2431& S3\cr
  & 64 & 100 & 0.009 & 0.39 & 2.13 &225 & 15 & 3375 & all\cr
    & 128 & 400 & 0.009 & 0.8 & 2.67 &282 & 28 & 7896 & all\cr
\hline
0.7 & 32 & 50 &0.097 & 0.29 & 5.81 & 61 & 11 & 671& all\cr
   &  32 & 100 &0.026 & 0.52 & 10.44 & 41 & 14 & 574 & S3\cr
  & 64 & 200 & 0.097 & 0.38 & 7.08 & 74 & 21 & 1554 & all\cr
    & 128 & 400 & 0.036 & 0.5 & 4.82 &135 & 29 & 3915 & all\cr
\hline
1 &32& 15 & 0.058 & 0.27 & 0.99 & 18& 7& 126 &  all \cr
  &32& 20 & 0.078 & 0.31 & 1.09 & 15& 8& 120 &  S3\cr
  &64 & 60 &0.058 & 0.35 & 1.05 & 19 & 13 & 247 &   all\cr
  & 128 & 240 & 0.058 & 0.64 & 1.05 & 19 & 23 & 437 & all\cr
\hline
\end{tabular}
\end{center}
\label{default1}
\end{table}

\begin{figure}[!h]
	\centering      	
	\includegraphics[width=7.0cm]{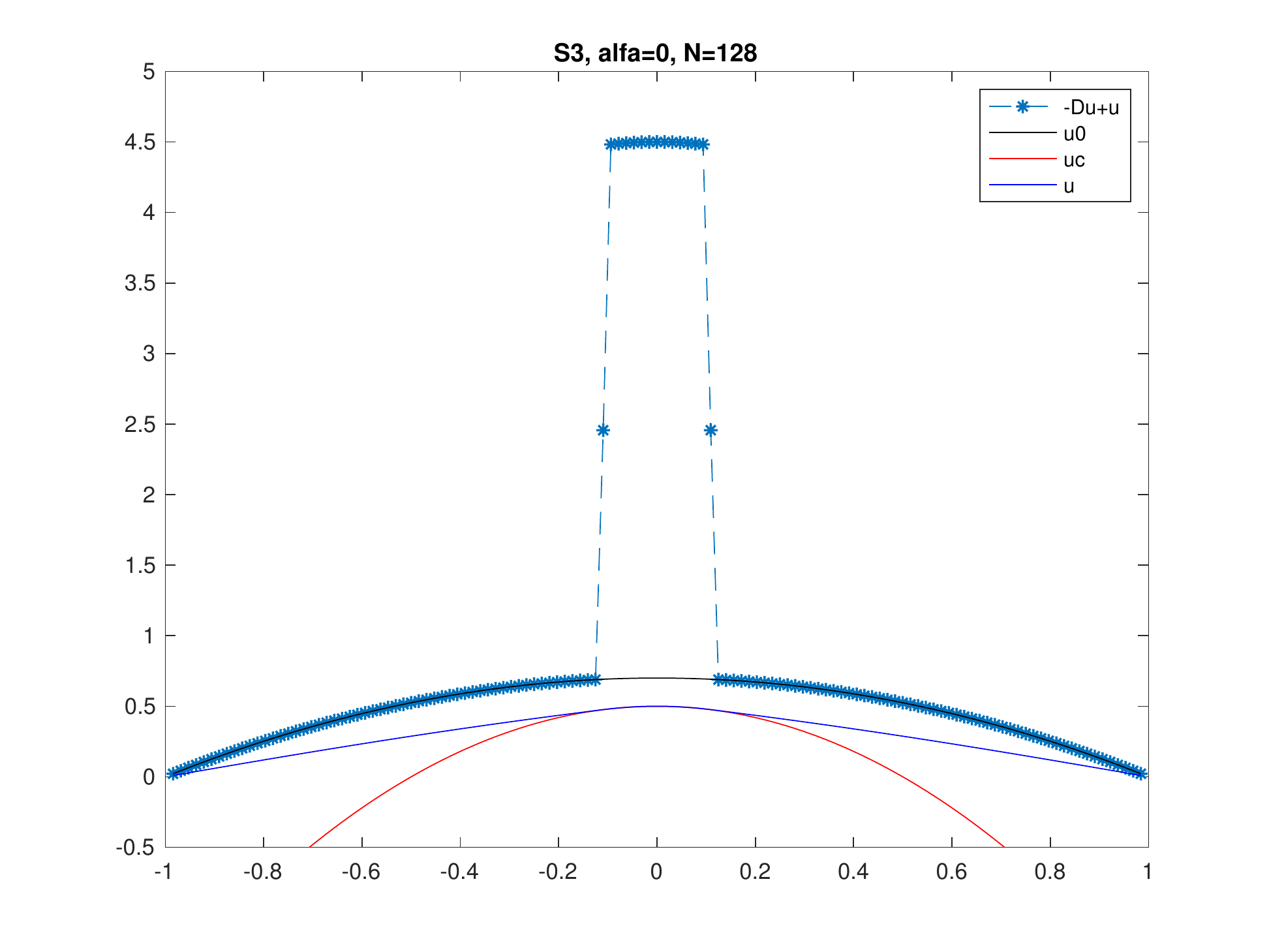}
	%\vskip-1.5cm
\caption{Solution for $\alpha=0$ and 128 nodes, scheme S3.}
	\label{F2}
\end{figure}

\begin{itemize}
\item The stationary solution is the same for any  $\alpha$, as stated by Theorem \ref{main}, and corresponds to the one of the stationary problem (\ref{BS}). Different is only the speed of convergence. The detected  right extremum of the symmetric contact set $C$ on the used meshes is 0.125  (the continuous value should be approximately $0.132$). 
\item Schemes S1 and S2 have a very similar behaviour; in both cases for  large values of $\gamma_\alpha$ the  contact set can be overestimated, a problem not present for  S3, due to the implicit nature of the quasi-Newton scheme of \cite{BS}. The semi implicit schemes S1 and S2 pay for the delay with which the contact information with the obstacle is achieved, allowing an uncorrect  evolution of the solution. To avoid that, strong $\tau$-step restrictions are necessary:  experiments show that all the schemes work correctly for example if $\tau^\alpha<0.1$, a bound which becomes particularly heavy  when $\alpha$ is small.  In the Table we reported approximately for any number of nodes the largest values of $\gamma_\alpha$ for which all the three schemes give the same correct solution.

\item Computational costs: the previous remark suggests that S3 is the more reliable and even the less expensive of the three schemes, allowing larger time steps and hence less iterations. Anyway, any single time iteration of S3 is much more expensive, many Picard iterations (growing with $N$)  with respect to a single linear system necessary to be solved for S1 and S2. Then,  comparing the computational cost of the schemes, even these two schemes reveal competitive in terms of the total number of linear system solved.

\item For a fixed  $\alpha$ the full contact time grows with the number of nodes, and does not seem to depend from $\gamma_\alpha$. On the contrary the stabilization time grows with $\gamma_\alpha$ but decreases with the number of nodes, since a finer mesh reduces the error at the boundary of the contact set.
\item All the schemes correctly work also for the case $\alpha=1$: it easy to see that  in that case both the used approximations of the Caputo derivative reduce to the standard incremental ratio in time.
\item The case $\alpha=0$ is a sort of control test: since in that case $\partial^\alpha_t u=u-u^0$, in absence of an obstacle the equation (\ref{model1}) would reduce to the stationary equation 
\begin{equation}\label{asint0}
-\Delta \o+\o=u^0 .
\end{equation}
In presence of an obstacle we expect  the discrete solution to satisfy  (\ref{asint0})  only outside the contact set, and from the first iteration. It is in fact clear from  (\ref{coefs1}) and (\ref{coefs4}) that since $H(u^0-\psi)=1$ and $\gamma_\alpha=1/h^2$,  the first iteration of all the schemes  becomes
$$ B^0 u^1= \left(I+\frac 1 {h^2}A\right)u^1 = u^0\ ,$$
which is essentially the discrete version of (\ref{asint0}). If $u^1$  goes over the obstacle, such identity will be satisfied only where  $u^1>\psi$. In Figure~\ref{F2} it can be seen what happens in our example with $N=128$ nodes and scheme S3: the solution $u$ is plotted in blue, the obstacle in red,  the initial datum $u^0$ in black and the quantity  $-\Delta u+u$ through asterisks:  the last two quantities coincide  in the detachment set, with a  natural discontinuity at the boundary of such a region.  Since the time step $\tau$ has no effect on the solution,  the  semi implicit schemes $S1$ and $S2$ for this example always overestimate the contact set.
\end{itemize}

%\vskip-2cm

\begin{figure}[!h]
	\centering      	
	\includegraphics[width=7.0cm]{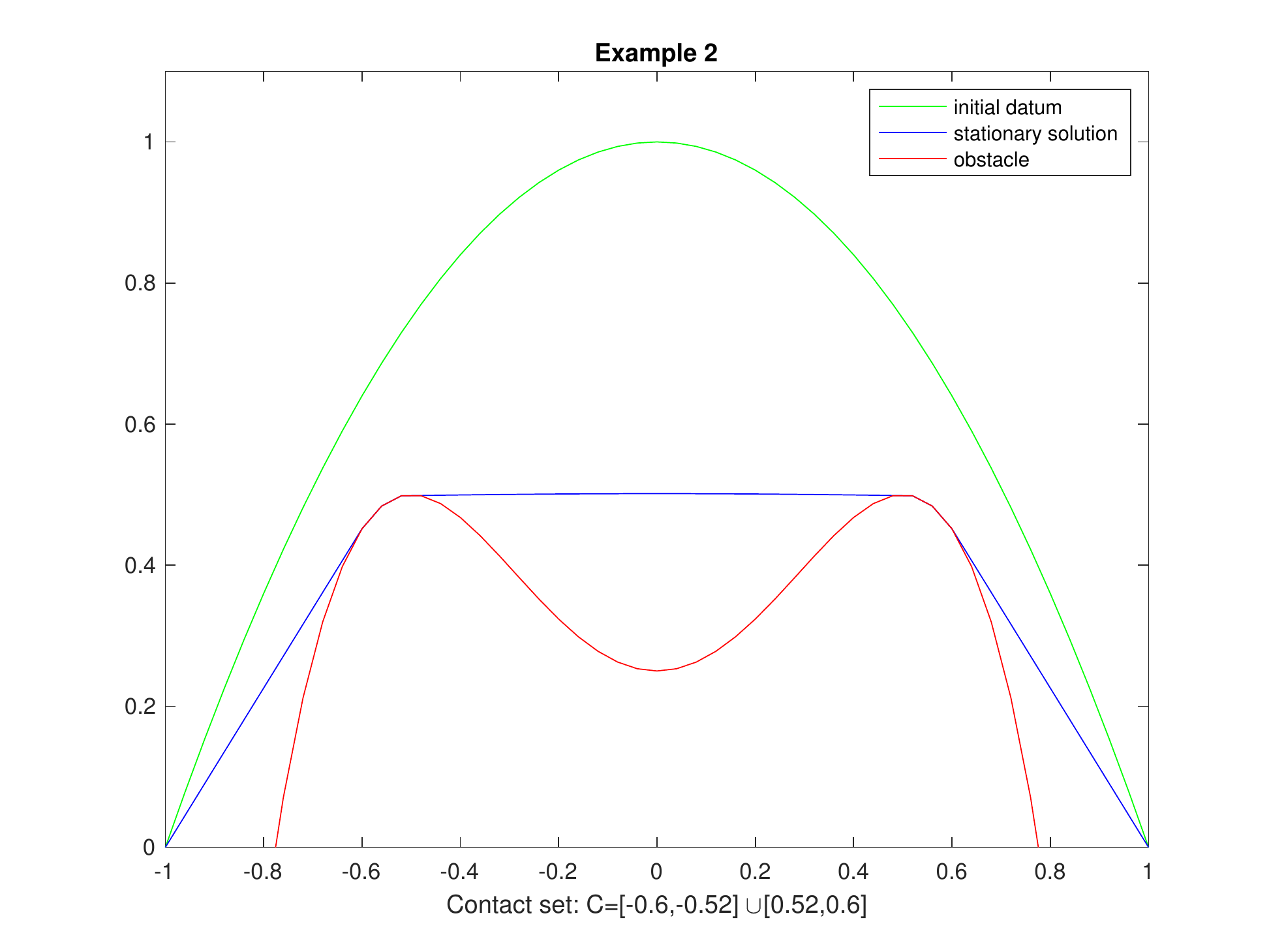}
	%\vskip-1.5cm
\caption{Data of Example 2.}
	\label{Ex2}
\end{figure}

\medskip\noindent {\bf Example 2.} $\Omega=(-1,1)$, $u^0(x)=1-x^2$, $\psi(x)=0.5-(2x^2-0.5)^2$ (see Fig.\ref{Ex2}).

\begin{figure}[ht!]
	\centering      	
	\includegraphics[width=7.5cm]{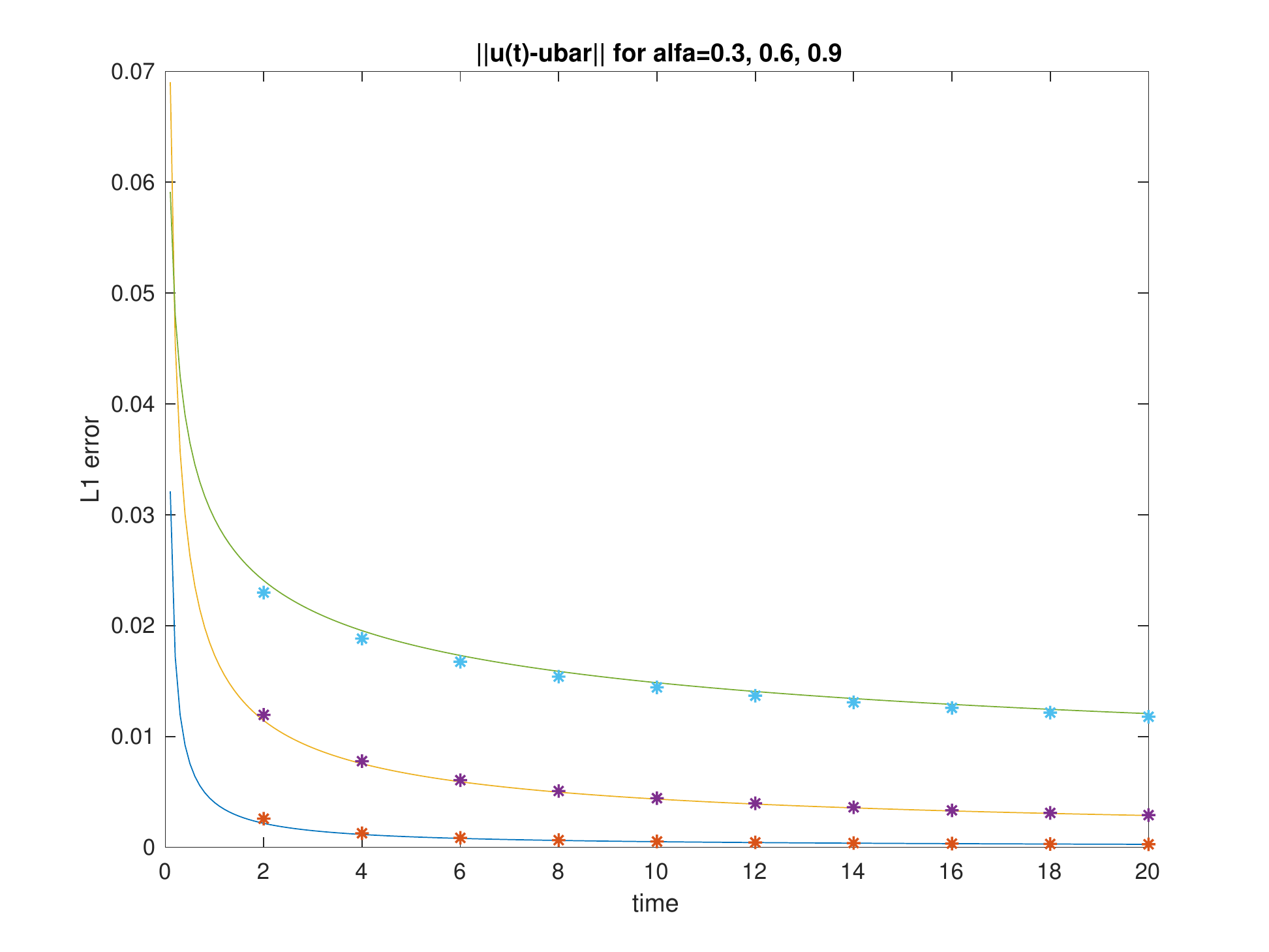}
	%\vskip-1.5cm
\caption{Example 2: polynomial speed of stabilization for $\alpha=0.3, 0.6$ and $0.9$.}
	\label{decay}
\end{figure}

\medskip On this example we tested numerically the estimate (\ref{stima1}) of Theorem~\ref{main}, using scheme S3. The $L^1$ norm of the error  at time $t^m$ on the given mesh was approximated by a natural quadrature formula, that is
$$\left\Vert u(t^m) - \overline{u}\right\Vert_{L^1\left(\Omega\right)}\simeq h\sum_i |u^m_i-\bar{u}_i| \ $$
(the vector $\bar{u}$ on the mesh was computed in advance with a sufficiently high precision).
Such an error  was computed for different values of $\alpha$  and of the time $t^m$. In Table \ref{default2} we reported  (in the third column) the discrete $L^1$ error  at $T=10$ with $N=32$ nodes and the same $\gamma_\alpha=50$ for any $\alpha$.  In the fourth column  the corresponding values of quantity $J(T)=\frac1{C}\frac{T^{-\alpha}}{\Gamma(1-\alpha)}$ of (\ref{jei}) are shown; for the constant $C$ we adopted a computed  estimate ($C=26$). It is evident that the two quantities decay at the same rate. In the second column it is also possible to see the number of instant times $M$ needed for each $\alpha$ to keep the same $\gamma_\alpha$, and the consequent increasing complexity of the computation.
In order to confirm  the order of decay of the error, polynomial for $\alpha\in(0,1)$ and exponential (up to the machine precision) for $\alpha=1$, we illustrate in Figure \ref{decay} such behavior plotting the continuous curves $J(t)$ in the time interval $(0,20)$ for $\alpha=0.3, 0.6$ and $0.9$, and the computed errors on ten different discrete times (marked by asterisks). These  values correctly lie on the curves.

%\begin{table}[ht]
%\caption{Example 2: $N=32$, $\gamma_\alpha=50$, scheme $S3$.}
%\begin{center}
%\begin{tabular}{|c|c|c|c|c|}
%\hline
%$\alpha$ & $M$ &  ${T^{-\alpha}}/{\Gamma(1-\alpha)}$ & error at $T=10$ & $~C$  \cr
%\hline
%1 & 51 &  0 & $1.72\ 10^{-15}$ &    \cr
%\hline
%0.9 & 61 & 0.0132 &  $6.84\ 10^{-6}$ &1900 \cr
%\hline
%0.8 & 77 &  0.0345 & $1.77\ 10^{-5}$ &1950 \cr
%\hline
%0.7 & 103 &  0.0667 & $3.40\ 10^{-5}$ &1960 \cr
%\hline
%0.6 & 152 &  0.1132 & $5.77\ 10^{-5}$ &1960 \cr
%\hline
%0.5 & 262 &  0.1784 & $9.10\ 10^{-5}$ &1960 \cr
%\hline
%0.4 & 593 &  0.2673 & $1.36\ 10^{-4}$ &1965 \cr
%\hline
%0.3 & 2313 &  0.3861 & $1.98\ 10^{-4}$ &1950 \cr
%\hline
%0.2 & 35184 &  0.5420 & $2.76\ 10^{-4}$ (est.)&1964 \cr
%\hline
%0.1 & 123000000 &  0.7433 & $3.79\ 10^{-4}$ (est.)&1961 \cr
%\hline
%\end{tabular}
%\end{center}
%\label{default2}
%\end{table}

\begin{table}[!h]
\caption{Example 2: $N=32$, $\gamma_\alpha=50$, scheme $S3$.}
\begin{center}
\begin{tabular}{|c|c|c|c|}
\hline
$\alpha$ & $M$ &   $L^1$ error at $T=10$ & ${T^{-\alpha}}/{(C\Gamma(1-\alpha))}$   \cr
\hline
1 & 51 & $8.07\ 10^{-6}$ & $\simeq 0$  \cr
\hline
0.9 & 61 & $5.27\ 10^{-4}$ &  $5.09\ 10^{-4}$  \cr
\hline
0.8 & 77 & $1.37\ 10^{-3}$ & $1.32\ 10^{-3}$  \cr
\hline
0.7 & 103 &$2.62\ 10^{-3}$& $2.56\ 10^{-3}$  \cr
\hline
0.6 & 152 & $4.41\ 10^{-3}$& $4.35\ 10^{-3}$  \cr
\hline
0.5 & 262 &  $6.87\ 10^{-3}$& $6.86\ 10^{-3}$  \cr
\hline
0.4 & 593 & $1.01\ 10^{-2}$ & $1.02\ 10^{-2}$  \cr
\hline
0.3 & 2313 &  $1.44\ 10^{-2}$  & $1.48\ 10^{-2}$ \cr
\hline
0.2 & 35184 &  $1.98\ 10^{-2}$ & $2.08\ 10^{-2}$  \cr
\hline
0.1 & 123000000 &  $2.7\ 10^{-2}$ (estimate) & $2.85\ 10^{-2}$   \cr
\hline
\end{tabular}
\end{center}
\label{default2}
\end{table}

%\newpage


\begin{thebibliography}{23}
\bibitem{ACF} 
Alberini C., Capitanelli R. and Finzi Vita S., 
{\em A numerical study of an Heaviside function driven degenerate diffusion equation}, 
preprint 2020, submitted.

\bibitem{Baz2000} 
Bazhlekova E.G.,
{\em Subordination principle for fractional evolution equations}, 
Fractional Calculus and Applied Analysis. {\bf 3}, No. 3, (2000), 213--230.

\bibitem{BS} 
Brugnano L. and Sestini A., 
{\em Iterative solution of piecewise linear systems for the numerical solution of obstacle problems}, 
J. Num. Anal. Ind. Appl. Math. (JNAIAM), {\bf 6}, No. 3-4, (2011), 67--82.

\bibitem{CapDov}
Capitanelli R. and D'Ovidio M., 
\emph{Fractional equations via convergence of forms}, 
Fractional Calculus and Applied Analysis, {\bf 22}, (2019), 844--870.

\bibitem{Dieth} 
Diethelm K., 
\emph{The analysis of fractional differential equations. An application-oriented exposition using differential operators of Caputo type}, 
Lecture Notes in Mathematics, 2004. Springer-Verlag, Berlin, (2010).

\bibitem{Dov}
D'Ovidio M., 
\emph{On the fractional counterpart of the higher-order equations},
Statistics and Probability Letters, {\bf 81}, (2011), 1929--1939.

\bibitem{giga} 
Giga Y., Liu Q. and Mitake H., 
{\em On a discrete scheme for time fractional fully nonlinear evolution equations}, 
Asymptotic Analysis, {\bf 120}, No. 1-2, (2020), 151--162.   %arXiv:1902.08863v1 [math.AP] (2019)

\bibitem{GKMR-book}
Gorenflo R., Kilbas A.A., Mainardi F. and  Rogosin S.V., 
\emph{Mittag-Leffler Functions, Related Topics and Applications},
Springer Monographs in Mathematics, Springer-Verlag Berlin Heidelberg (2014).

\bibitem{JLZ} 
Jin B, Lazarov R. and Zhou Z., 
{\em Numerical methods for time-fractional evolution equations with nonsmooth data: a concise overview},  
Comput. Methods Appl. Mech. Engrg., {\bf 346}, (2019), 332--358. 

\bibitem{Koc89} 
Kochubei  A.N.,
{\em The Cauchy problem for evolution equations of fractional order},
Differential Equations, {\bf 25}, (1989), 967--974.


\end{thebibliography}
\end{document}